\documentclass[11pt,a4paper]{article}

\usepackage{latexsym,amssymb,amsfonts,amsmath,epsfig,tabularx,stmaryrd,amsthm}
\usepackage{mathtools}
\usepackage[usenames]{xcolor}
\usepackage{graphicx}

\newtheorem{theorem}{Theorem}
{\theoremstyle{definition}

\newtheorem{remark}{Remark}
}
\newtheorem{lemma}[theorem]{Lemma}
\newtheorem{corollary}[theorem]{Corollary}

\usepackage{graphicx}       

\usepackage[latin1]{inputenc}
\usepackage[T1]{fontenc}
\usepackage{microtype} 

\usepackage{cite}
\usepackage{array,colortbl}
\DeclareFontFamily{U}{mathx}{\hyphenchar\font45}
\DeclareFontShape{U}{mathx}{m}{n}{
      <5> <6> <7> <8> <9> <10>
      <10.95> <12> <14.4> <17.28> <20.74> <24.88>
      mathx10
      }{}
\DeclareSymbolFont{mathx}{U}{mathx}{m}{n}
\DeclareFontSubstitution{U}{mathx}{m}{n}
\DeclareMathAccent{\widecheck}{0}{mathx}{"71}

\def\citep#1#2{\cite[{#1}]{#2}}

\usepackage{algorithm}
\usepackage{algorithmic}

\newcommand{\bbE}{{\mathbb{E}}}

\newcommand{\bbN}{{\mathbb{N}}}

\newcommand{\bbP}{{\mathbb{P}}}

\newcommand{\bbR}{{\mathbb{R}}}

\newcommand{\N}{{\mathbb{N}}} 
\newcommand{\R}{{\mathbb{R}}} 

\DeclareSymbolFont{bbold}{U}{bbold}{m}{n}
\DeclareSymbolFontAlphabet{\mathbbold}{bbold}
\newcommand{\ind}{{\mathbbold{1}}}

\newcommand{\bse}{{\boldsymbol{e}}}

\newcommand{\bsy}{{\boldsymbol{y}}}
\newcommand{\bsz}{{\boldsymbol{z}}}
\newcommand{\bsA}{{\boldsymbol{A}}}

\newcommand{\bsY}{{\boldsymbol{Y}}}

\newcommand{\bszero}{{\boldsymbol{0}}} 
\newcommand{\bsone}{{\boldsymbol{1}}}  
\newcommand{\bsalpha}{{\boldsymbol{\alpha}}}

\newcommand{\bsgamma}{{\boldsymbol{\gamma}}}

\newcommand{\bseta}{{\boldsymbol{\eta}}}

\newcommand{\bsmu}{{\boldsymbol{\mu}}}
\newcommand{\bsnu}{{\boldsymbol{\nu}}}

\newcommand{\bsrho}{{\boldsymbol{\rho}}}

\newcommand{\bstau}{{\boldsymbol{\tau}}}

\newcommand{\bspsi}{{\boldsymbol{\psi}}}

\newcommand{\bsDelta}{{\boldsymbol{\Delta}}}

\newcommand{\bsPhi}{{\boldsymbol{\Phi}}}

\newcommand{\calH}{{\mathcal{H}}}

\newcommand{\calO}{{\mathcal{O}}}

\newcommand{\setu}{{\mathfrak{u}}}

\newcommand{\rd}{\,\mathrm{d}} 

\graphicspath{{.}{./Figures/}{../CDF_approx/Figures/}}

 {\theoremstyle{plain}\newtheorem{GKSSassumption}{Assumption}} 

\newcommand{\GKSSind}{\mathop{\rm ind}}
\newcommand{\GKSSHe}{\mathrm{He}}

\title{Theory and construction of Quasi-Monte Carlo rules for Asian option pricing and density
  estimation}
\author{Alexander D.~Gilbert, Frances Y.~Kuo, Ian H.~Sloan and Abirami Srikumar%
   \footnote{School of Mathematics and Statistics, UNSW Sydney, Sydney NSW 2052, Australia \newline
   Emails: alexander.gilbert@unsw.edu.au, f.kuo@unsw.edu.au, i.sloan@unsw.edu.au,
   a.srikumar@student.unsw.edu.au}
   }
\date{\today}

\begin{document}

\maketitle

\abstract{In this paper we propose and analyse a method for estimating
three quantities related to an Asian option: the fair price, the
cumulative distribution function, and the probability density.  The method
involves preintegration with respect to one well chosen integration
variable to obtain a smooth function of the remaining variables, followed
by the application of a tailored lattice Quasi-Monte Carlo rule to
integrate over the remaining variables.}

\section{Introduction}

This paper proposes and analyses a preintegration
method for estimating \emph{the fair price of an Asian option, and the
associated distribution function and density function}.  The method
approximates the corresponding multiple integrals by first preintegrating
with respect to one well chosen variable, resulting in a smooth function
of the other variables, and then integrates over the remaining variables
using a carefully constructed Quasi-Monte Carlo method of lattice rule
type.

An Asian option is an option contract where the payoff is determined by
the arithmetic average of the price of the underlying asset over the time
period. The payoff for an Asian \emph{put} option is
\begin{equation*}
\mathrm{payoff} = \max(K - X,0),
\end{equation*}
where $K$ is the fixed strike price and $X$ is a real-valued random
variable for the average price. Similarly, the \emph{call} option is
$\max(X - K, 0)$. Due to this symmetry, 
it suffices to consider only one of the two. Under the Black--Scholes model, the (discretised)
average price is given by a smooth function $X = \phi(Y_0, Y_1, \ldots,
Y_d)$, where $Y_\ell \sim \mathrm{N}(0, 1)$ i.i.d.
The definition of $\phi$ depends on how the path of the 
underlying stock price is generated and so is not unique.
This choice will affect how preintegration can be applied and
will be discussed in Section~\ref{GKSS:sec:option}.

In this paper, we are interested in computing three quantities  related to
the Asian option, namely the \emph{value} or \emph{fair price} of the
option; the \emph{probability} that the average price is below a certain
value;  and the \emph{probability density} of the random variable $X$.
Each of the problems can be formulated as an expected value of some
non-smooth function of the average price, or equivalently, as a
high-dimensional integral, as follows:
\begin{itemize}
\item The value of the Asian put option is the discounted expected
    value of the payoff
\begin{equation}
\label{GKSS:eq:price}
\text{value} = e^{-R T}\bbE[ \max(K - X, 0)]
= e^{-R T}\int_{\R^{d+1}} \max(K - \phi(\bsy), 0) \bsrho(\bsy) \rd \bsy,
\end{equation}
where $R$ is the risk-free interest rate, $T$ is the expiration time,
and $\bsrho(\bsy)$ is the product standard normal density.  Recall
that $K$ is the strike price and $\phi(\bsy)$ is the average price.

\item The probability that the average is below $x \in \R$ is known as
    the \emph{cumulative distribution function (cdf)} of $X$ at $x$,
    which is denoted by $F : \R \to [0, 1]$. It can be formulated as
    the expected value of an indicator function of the average price
\begin{equation}
\label{GKSS:eq:cdf}
F(x) = \bbP[X \leq x]
= \bbE[\GKSSind(x - X)] =
\int_{\R^{d+1}} \GKSSind(x - \phi(\bsy)) \bsrho(\bsy) \rd \bsy,
\end{equation}
where $\GKSSind$ is the indicator function defined by $\GKSSind(z) \coloneqq
1$ if $z \geq 0$ and~$0$ otherwise.

\item The \emph{probability density function (pdf)}, denoted by $f :
    \R \to [0, \infty)$, can be formulated as the expected value of a
    Dirac $\delta$ of the average price
\begin{equation}
\label{GKSS:eq:pdf}
f(x) = \bbE[ \delta( x - X)]
= \int_{\R^{d+1}} \delta(x - \phi(\bsy)) \bsrho(\bsy) \rd \bsy,
\end{equation}
where $\delta$ is the distribution satisfying $\delta(z)= 0$ for all
$z\ne 0$ and $\int_{-\infty}^\infty \delta(z) w(z) \rd z = w(0)$ for
all sufficiently smooth function $w$.
\end{itemize}
Each integral \eqref{GKSS:eq:price}--\eqref{GKSS:eq:pdf} can be
written as $\int_{\R^{d+1}} g(\bsy) \bsrho(\bsy) \mathrm{d} \bsy$, with the integrand given by
\begin{align}\label{GKSS:eq:int}
g(\bsy) =
\begin{cases}
 \max(K - \phi(\bsy), 0) = (K- \phi(\bsy))\,\GKSSind(K - \phi(\bsy)) & \mbox{for  option value}, \\
 \GKSSind(x - \phi(\bsy)) & \mbox{for cdf $F(x)$}, \\
 \delta(x - \phi(\bsy)) & \mbox{for pdf $f(x)$}.
 \end{cases}
\end{align}

Preintegration is a method of smoothing out a function with a simple
discontinuity or kink (i.e. discontinuity in the first derivative) by
integrating first with respect to a single specially chosen variable. We
denote this special variable by $y_0$ and the remaining $d$ variables by
$\bsy$, so that the preintegration step is to evaluate
\begin{equation}
\label{GKSS:eq:P_0}
P_0 g(\bsy) \coloneqq \int_{-\infty}^\infty g(y_0, \bsy) \rho(y_0) \rd y_0
 \qquad\mbox{for } \bsy\in\R^d.
\end{equation}
Each desired integral is then a $d$-dimensional integral $\int_{\R^d} P_0
g(\bsy) \rho(\bsy)\rd\bsy$, in which $P_0 g$ is smooth under suitable
conditions (see Section~\ref{GKSS:sec:preint}), of which the most
crucial is that $\phi(y_0, \bsy)$ be an increasing (or decreasing) function
of $y_0$ for each $\bsy \in \R^d$.
In practice, the preintegration variable $y_0$ should 
foremost be chosen such that these conditions are satisfied,
then secondary considerations such as variance reduction
can also be taken into account.

Preintegration can be considered as a special case of
\emph{conditional expectation}, a method
often used to reduce the variance of multivariate integrands,
see, e.g., \cite{GKSS:ACN13a,GKSS:Glasserman,GKSS:GlaSta01,GKSS:LEcLem00,GKSS:LEcPuchBAb22,GKSS:WWH17}. 
As a method of smoothing the integrand it was
presented, with supporting theory, in \cite{GKSS:GKLS18}.  In turn, that paper
built upon earlier results
\cite{GKSS:GKS10, GKSS:GKS17note} concerning the smoothing effect of integration.
Recently, \cite{GKSS:Liu22,GKSS:LiuOwen23} studied how 
a change of variables can be performed such that preintegration
results in good variance reduction.

\section{Background}

\subsection{Asian options}
\label{GKSS:sec:option}

In the Black--Scholes model the price $S_t$ satisfies the \emph{stochastic
differential equation} (SDE)  $\rd S_t = R S_t \rd t + \sigma S_t \rd W_t
$ for $t \geq 0$ under the risk neutral measure, where $R > 0$ 
is the risk-free interest rate and
$\sigma > 0$ is the volatility, both of which are assumed to be
deterministic and constant in time, and $W_t$ is a standard Brownian
motion. The solution $S_t$ is given explicitly by
\begin{equation*}
S_t = S_0 \exp\big((R - \tfrac{1}{2}\sigma^2) t + \sigma W_t\big),
\end{equation*}
which in practice is simulated by discretising the time interval $[0, T]$
using $d + 1$ uniform timesteps, and then generating a $(d+1)$-dimensional
discrete Brownian motion given by a multivariate normal with covariance
matrix $\Sigma$. Then the average price over $[0, T]$ can be computed
by averaging the price at the discrete timesteps.

After factorising the covariance matrix as $\Sigma = AA^\top$ and making a
change of variables, the time-discretised average of the underlying asset
price $X = \phi(\bsY)$ is given by a function $\phi:\R^{d+1} \to \R$, of
$(d + 1)$ i.i.d.\ standard normals, $Y_\ell \sim \mathrm{N}(0, 1)$, with
\begin{equation}
\label{GKSS:eq:phi}
  \phi(\bsy)
= \frac{1}{d+1} \sum_{k = 0}^d S_0
\exp \bigg( \big(R - \tfrac{1}{2}\sigma^2\big)
\frac{(k+1)T}{d+1} + \sigma \bsA_k \bsy\bigg).
\end{equation}
where $\bsy \coloneqq (y_0, y_1, \ldots, y_d)$ and $\bsA_k \coloneqq
(A_{k, i})_{i = 0}^d$.

The factorisation of the covariance matrix, $\Sigma = AA^\top$, is not
unique, and different factorisations will affect the structure of the
problem, in turn affecting the performance of different quadrature
methods. The three most popular factorisation methods are Cholesky
factorisation, Brownian bridge, and principal components analysis (PCA),
see \cite{GKSS:Glasserman}. In this paper, we focus on the PCA factorisation,
which is based on the spectral decomposition of $\Sigma$
and orders the variables according to the eigenvalues of $\Sigma$ in decreasing value. For $\Sigma$
corresponding to $d + 1$ uniform timesteps, the entries of PCA factors are,
for $k,i = 0,1,\ldots,d$,
\begin{equation}
\label{GKSS:eq:A_pca}
 A_{k,i} = \tau_d\,\frac{\sin(2(k+1) \chi_i)}{\sin(\chi_i)},
 \; \tau_d \coloneqq \sqrt{\frac{T}{(d+1)(2d + 3)}},
 \; \chi_i \coloneqq \frac{\pi(2i + 1)}{2(2d + 3)} \in (0,\tfrac{\pi}{2}).
\end{equation}

In Section~\ref{GKSS:sec:qmc-preint}, we show that for the PCA factorisation 
\eqref{GKSS:eq:A_pca}, if $y_0$ is chosen as the preintegration variable then
$\phi$ \eqref{GKSS:eq:phi} satisfies the necessary monotonicity condition (cf.~Assumption~\ref{GKSS:asm:phi})
and is in fact the only variable for which $\phi$ is monotone.
Furthermore, since $y_0$ corresponds to the largest eigenvalue of $\Sigma$, we
can expect that it also offers the greatest
reduction in variance.

\subsection{Function spaces}

We define the needed function spaces in $(d + 1)$ dimensions; the
analogous spaces in $d$ dimensions are obtained simply by omitting $y_0$.
For $\bsnu \in \N^{d + 1}$, we write the $\bsnu$th mixed derivative as
$ D^\bsnu \coloneqq \prod_{i = 0}^d \partial^{\nu_i}/\partial y_i^{\nu_i},
$ and define the space of functions with continuous mixed derivatives up
to order $\bsnu$ by
$ C^\bsnu(\R^{d+1}) \,\coloneqq\, \{ \varphi \in C(\R^{d+1}):  D^\bseta
\varphi \in C(\R^{d+1}) \text{ for all } \bseta \leq \bsnu\}. $
Also, for a collection of weight parameters $\gamma_\setu > 0$ for $\setu
\subseteq \{0, 1, \ldots, d\}$ and weight functions $\psi_i : \R \to \R_+$
for $i = 0, 1, \ldots, d$ we define the weighted \emph{Sobolev space of
dominating mixed smoothness} of order $\bsnu$, denoted by $\calH_{d +
1}^\bsnu$, to be the space of locally integrable functions on $\R^{d + 1}$
such that the norm
\begin{equation*}
\|h\|_{\calH^\bsnu_{d + 1}}^2 \,\coloneqq\,
\sum_{\bseta \leq \bsnu} \frac{1}{\gamma_\bseta}
\int_{\R^{d + 1}} |D^\bseta h(y_0, \bsy) |^2 \bspsi_\bseta(\bsy_\bseta)
\bsrho_{-\bseta}(\bsy_{-\bseta})
\rd y_0 \rd \bsy
\end{equation*}
is finite. Here we have used the notation $\gamma_\bseta \coloneqq
\gamma_{\mathrm{supp}(\bseta)}= (y_i : \eta_i \neq 0)$,
$\gamma_{-\bseta}\coloneqq (y_i : \eta_i = 0)$, as well as the products
$\bspsi_{\bseta}(\bsy_\bseta) \coloneqq \prod_{i = 0,\, \eta_i \ne 0}^d
\psi_i(y_i)$ and $\bsrho_{-\bseta}(\bsy_{-\bseta}) \coloneqq \prod_{i =
0,\, \eta_i = 0}^d \rho(y_i)$. Throughout we assume that the weight
functions satisfy
\begin{equation}
\label{GKSS:eq:psi-cond}
 \int_{-\infty}^\infty \psi_i(y)\rd y = 1
 \quad\mbox{and}\quad
\int_{-\infty}^\infty \frac{\Phi(y)(1 - \Phi(y))}{\psi_i(y)} \rd y < \infty,
\end{equation}
where $\Phi$ is the cdf corresponding to the density $\rho$.

\subsection{Quasi-Monte Carlo methods}
\label{GKSS:sec:qmc}

A Quasi-Monte Carlo (QMC) approximation of an integral over $\R^d$ is an
equal-weight quadrature rule
\begin{equation}
\label{GKSS:eq:qmc}
\int_{\R^d} h(\bsy) \bsrho(\bsy) \, \rd \bsy
\,\approx\,
Q_{d, N} (h) \,\coloneqq\, \frac{1}{N} \sum_{k = 0}^{N - 1} h(\bsPhi^{-1}(\bstau_k))\,,
\end{equation}
where $N$ is the number of function evaluations and $\{\bstau_k\}_{k =
0}^{N - 1} \subset~[0, 1]^d$ are prescribed quadrature points, which
are mapped to $\R^d$ by applying the inverse cdf $\Phi^{-1}$
corresponding to $\rho$ componentwise.
\emph{Randomly shifted lattice rules} are a simple yet powerful class of
QMC methods where the whole pointset is generated by a single integer
\emph{generating vector} $\bsz \in \N^d$. The points for a randomly
shifted (rank-1) lattice rule on the unit cube are
\begin{equation}
\label{GKSS:eq:lattice}
\bstau_k \,=\, \bigg\{\frac{k\bsz}{N} + \bsDelta \bigg\} \in [0, 1]^d, \quad \text{for } k = 0, 1, \ldots, N - 1,
\end{equation}
where $\{ \cdot \}$ denotes taking the fractional part of each component
and $\bsDelta \sim \mathrm{Uni}[0, 1)^d$ is a single uniformly distributed
\emph{random shift}.

Good generating vectors and in turn, good lattice rules can be constructed
using the \emph{component-by-component} (CBC) construction. In \cite{GKSS:NK14}
it was shown that CBC-generated lattice rule approximations in $\R^d$
obtain the optimal rate of convergence of close to $1/N$ for the
error. For $h \in \calH_{d}^\bsone$, the root-mean-square error (RMSE) of a
CBC generated randomly shifted lattice rule approximation
\eqref{GKSS:eq:qmc}--\eqref{GKSS:eq:lattice} is bounded by 
\begin{align}
\label{GKSS:eq:cbc-err}
&\text{RMSE} = \sqrt{\bbE_\bsDelta\, \bigg|\int_{\R^d} h(\bsy) \bsrho(\bsy) \rd \bsy
- Q_{d, N}(h)\bigg|^2 } \nonumber\\
&\leq [\varphi_{\mathrm{tot}}(N)]^{-\frac{1}{2\lambda}}\,
 \bigg(\sum_{\bszero \neq \bseta \in \{0, 1\}^d}
\gamma_\bseta^\lambda \, [2 C_2 \zeta(2\omega\lambda)]^{|\bseta|}
\bigg)^{\frac{1}{2\lambda}}\,
\|h\|_{\calH_{d}^{\bsone}}
\quad\forall\, \lambda\in \big(\tfrac{1}{2\omega}, 1\big].
\end{align}%
Here $\zeta$ is the Riemann zeta function, $\varphi_{\rm tot}$ is the
Euler totient function, and $C_2 < \infty$, $\omega  = r_2 \in (1/2, 1)$ 
are as in \cite{GKSS:NK14}. Note that we follow \cite{GKSS:GKS22b} in
formulating the error in terms of $\calH_{d}^\bsone$ instead of the
\emph{unanchored ANOVA} space from \cite{GKSS:NK14}; see \cite{GKSS:GKS22a}, which
proves the required equivalence under the condition \eqref{GKSS:eq:psi-cond}.

\section{Smoothing by preintegration}
\label{GKSS:sec:preint}

In this section we review the theory from \cite{GKSS:GKS22b} on using smoothing
by preintegration for density estimation, which covers the cdf and pdf in
\eqref{GKSS:eq:int}. We also extend this theory to cover the option price
in the first case of \eqref{GKSS:eq:int}. To this end, we follow
\cite{GKSS:GKLS18,GKSS:GKS22b} and consider a generic function of the form
\begin{equation}
\label{GKSS:eq:g_ind}
g(y_0, \bsy) = \theta(y_0, \bsy) \GKSSind(x - \phi(y_0, \bsy))
\quad \text{or} \quad
g(y_0, \bsy) = 
\delta(x - \phi(y_0, \bsy)),
\end{equation}
where $x \in \R$, $\theta, \phi : \R^{d+1} \to \R$ and 
the components of $\bsy$ are i.i.d.~with marginal density~$\rho$. 
The option price fits the first framework with $\theta = x - \phi$ and $K = x$, while the
cdf fits with $\theta \equiv 1$. The pdf fits the second framework.

Throughout this section we will treat $\theta, \phi$, as generic functions
and $\rho$ as a generic density function. We will make assumptions about
them in order to state bounds on the norm of the preintegrated function
$P_0 g$. In the next section we will verify these assumptions for our
three specific cases in \eqref{GKSS:eq:int}.

\begin{GKSSassumption}
\label{GKSS:asm:phi} \textit{For $d\ge 1$ and $\bsnu\in\bbN_0^d$, let
$\theta, \phi: \R^{d+1} \to \R$ and $\rho\in \bbR\to\bbR$ satisfy}
\begin{enumerate}
\item \textit{$D^0\phi(y_0,\bsy) > 0$ for all $(y_0,\bsy) \in
    \R^{d+1}$;}
\item \textit{for each $\bsy \in \R^d$, $\phi(y_0, \bsy) \to \infty$
    as $y_0 \to \infty$;}
\item \textit{$\theta, \phi \in \calH_{d+1}^{(\nu_0, \bsnu)} \cap
    C^{(\nu_0, \bsnu)}(\R^{d+1})$, where $\nu_0 = |\bsnu| + 1$; and}
\item \textit{$\rho \in C^{\nu_0 - 1}(\R)$.}%
\end{enumerate}
\end{GKSSassumption}

The first two items imply that $\phi$ is strictly increasing with respect
to $y_0\in\bbR$ and also tends to $\infty$ as $y_0\to\infty$. For $x \in
[a,b]$ we define
\begin{equation}
\label{GKSS:eq:Ux}
 U_x \coloneqq \{\bsy \in \R^d : \phi(y_0, \bsy) = x
\; \text{for some } y_0 \in \R\}.
\end{equation}
Then the implicit function \cite[Theorem~3.1]{GKSS:GKS22b} states that if the
set
\begin{equation}
\label{GKSS:eq:V}
V \coloneqq \big\{(x,\bsy) \in (a, b) \times \R^d : \phi(y_0, \bsy)  = x
\text{ for some } y_0 \in \R \big\} \subset [a, b] \times \R^d
\end{equation}
is not empty, then there exists a unique function $\xi \in C^{(\nu_0,
\bsnu)}(\overline{V})$ satisfying
\begin{align} \label{GKSS:eq:xi}
  \phi(\xi(x,\bsy),\bsy) = x
  \qquad  \mbox{for all} \quad (x, \bsy) \in \overline{V}.
\end{align}
Furthermore, for $(x, \bsy) \in V$ the first-order derivatives of $\xi(x, \bsy)$ are given
by
\begin{align}
\label{GKSS:eq:Di-xi}
  D^i \xi(x,\bsy)
  &= - \frac{D^i\phi(\xi(x,\bsy),\bsy)}{D^0 \phi(\xi(x,\bsy),\bsy)}
  \qquad\mbox{for all } i=1,\ldots,d, \quad \text{and}
  \\
\nonumber
  \frac{\partial}{\partial x} \xi(x,\bsy)
  &= \frac{1}{D^0 \phi(\xi(x,\bsy),\bsy)}.
\end{align}

It then follows that the preintegrated function \eqref{GKSS:eq:P_0}
simplifies to
\begin{align} \label{GKSS:eq:P0g-simp}
&P_0g(\bsy)
\\\nonumber
&=\begin{cases}
 \displaystyle
 \int_{-\infty}^{\xi(x,\bsy)} \theta(y_0,\bsy) \rho_0(y_0) \rd y_0
 & \text{if } g(y_0,\bsy) = \theta(y_0, \bsy)\GKSSind(x - \phi(y_0,\bsy)),\\
 \Phi(\xi(x,\bsy))
 & \text{if } g(y_0,\bsy) = \GKSSind(x - \phi(y_0,\bsy)), \\
 \displaystyle \frac{\mathrm{\partial}}{\mathrm{\partial} x} \Phi(\xi(x,\bsy)) = \frac{\rho(\xi(x,\bsy))}{D^0\phi(\xi(x,\bsy),\bsy)}
 & \text{if } g(y_0,\bsy) = \delta(x - \phi(y_0,\bsy)),
\end{cases}
\end{align}
for $\bsy \in U_x$, and all cases are $0$ if $\bsy\in\bbR^d\setminus U_x$.

To establish the smoothness of $P_0g$, we need bounds involving the
typical functions that arise from taking mixed derivatives of $P_0g$ under
the chain rule. The following assumption generalises
\cite[Assumption~2]{GKSS:GKS22b} to include the extra function $\theta$ in
\eqref{GKSS:eq:g_ind}. (Taking $\theta \equiv 1$ and $\bsmu=\bszero$ recovers
\cite[Assumption~2]{GKSS:GKS22b}.)

\begin{GKSSassumption} \label{GKSS:asm:theta-phi} \textit{
Let $d\ge 1$, $\bsnu\in \bbN_0^d$, $[a,b]\subset\bbR$, and suppose that
$\theta$, $\phi$ and $\rho$ satisfy Assumption~\ref{GKSS:asm:phi}. Recall
the definitions of $U_x$, $V$ and $\xi$ in \eqref{GKSS:eq:Ux},
\eqref{GKSS:eq:V} and \eqref{GKSS:eq:xi}, respectively.
Given $q\in\bbN_0$ and $\bseta\le\bsnu$ satisfying $|\bseta|+q \le
|\bsnu|+1$, we consider functions $h_{q,\bseta}: \overline{V} \to \R$ of
the form
\begin{align} \label{GKSS:eq:h-form}
 \hspace{-0.5cm}
 \begin{cases}
 h_{q, \bseta}(x, \bsy) = h_{q,\bseta,(r,\bsalpha,\beta,\bsmu)}(x,\bsy) \\
 \hspace{1.4cm}
 \coloneqq
 \displaystyle\frac{(-1)^r D^\bsmu\theta(\xi(\bsy), \bsy)\,
       \rho^{(\beta)}(\xi(x,\bsy))\,
       \prod_{\ell=1}^{r} D^{\bsalpha_\ell} \phi(\xi(x,\bsy),\bsy)}
      {[D^0\phi(\xi(x,\bsy),\bsy)]^{r+q}},
\\[4mm]
 \mbox{with $r\in\bbN_0$, $\bsalpha=(\bsalpha_\ell)_{\ell=1}^r$,
 $\bsalpha_\ell\in \N_0^{d + 1}\!\setminus\!\{\bse_0,\bszero\}$, $\beta \in \N_0$, $\bsmu\in\N_0^{d+1}$ satisfying}
 \\[1mm]
 \hspace{0.5cm} r\le 2|\bseta|+q-1,\;
 \alpha_{\ell,0} \le |\bseta|+q, \;
 \beta \leq |\bseta|+q-1, \;
 |\bsmu| \le |\bseta|, \\
 \hspace{0.5cm} \bsmu + \beta\bse_0 + \displaystyle\sum_{\ell = 1}^{r} \bsalpha_\ell
 = (r+q-1, \bseta).
 \end{cases}
 \hspace{-0.5cm}
\end{align}
We assume that all such functions $h_{q,\bseta}$ satisfy
\begin{equation} \label{GKSS:eq:h-lim}
 \lim_{\bsy \to \partial U_x} h_{q,\bseta}(x, \bsy) = 0
 \quad \text{for all } x \in [a,b],
\end{equation}
and there is a constant $B_{q,\bseta}$ such that
\begin{equation} \label{GKSS:eq:h-int}
 \sup_{x \in [a, b]}
 \int_{U_x} |h_{q,\bseta}(x, \bsy)|^2 \,
 \bspsi_{\bseta}(\bsy_\bseta)\,\bsrho_{-\bseta}(\bsy_{-\bseta})\rd \bsy
 \le B_{q,\bseta} \,<\, \infty.
\end{equation}}%
\end{GKSSassumption}

The following lemma counts the number of terms that arise from the chain
rule and implicit differentiation when differentiating a function of the
form \eqref{GKSS:eq:h-form}. It will be key to obtaining an explicit bound
on the norm of $P_0g$. This lemma generalises \cite[Lemma~B.1]{GKSS:GKS22b} to
include the extra function $\theta$ in \eqref{GKSS:eq:h-form}.

\begin{lemma} \label{GKSS:lem:count}
Let $d\ge 1$, $\bsnu \in \bbN_0^d$, $[a,b]\subset\R$, and suppose that
$\theta,\phi$ and $\rho$ satisfy Assumption~\ref{GKSS:asm:phi}. Recall the
definitions of $U_x$, $V$ and $\xi$ in \eqref{GKSS:eq:Ux},
\eqref{GKSS:eq:V} and \eqref{GKSS:eq:xi}, respectively. For any
$q\in\bbN_0$ and $\bseta\le\bsnu$ satisfying $|\bseta|+q\le |\bsnu|+1$, we
consider functions $h_{q, \bseta} : V \to \R$ of the form
\eqref{GKSS:eq:h-form}. Then for any $i\in \{1:d\}$, we can write
\begin{equation}
\label{GKSS:eq:Dh-sum}
 D^i h_{q,\bseta}(x,\bsy)
 = \sum_{k=1}^{K_{q,\bseta}} h_{q,\bseta+\bse_i}^{[k]}(x,\bsy)
 \quad\mbox{with}\quad
 K_{q,\bseta} \le 8|\bseta| +6q-1,
\end{equation}
where each function $h_{q,\bseta+\bse_i}^{[k]}$ is of the form
\eqref{GKSS:eq:h-form} with $\bseta$ replaced by $\bseta+\bse_i$. For the
special case $\theta \equiv 1$ and $\bsmu = \bszero$ in
\eqref{GKSS:eq:h-form}, the number of functions reduces~to $K_{q,\bseta}
\le 8|\bseta| +6q-3$.
\end{lemma}

\begin{proof}
The case $\theta \equiv 1$ is given by \cite[Lemma~B.1]{GKSS:GKS22b}.
For the general case, we can write
\begin{equation}
\label{GKSS:eq:h=Dtheta-h}
h_{q, \bseta, (r, \bsalpha, \beta,\bsmu)}(x,\bsy)
= D^\bsmu \theta(\xi(x,\bsy), \bsy) \,
\widetilde{h}_{q, \bseta, (r, \bsalpha, \beta, \bsmu)}(x,\bsy),
\end{equation}
where we defined
\begin{equation}
\label{GKSS:eq:h-tilde}
\widetilde{h}_{q, \bseta, (r, \bsalpha, \beta, \bsmu)}(x,\bsy)
\coloneqq \frac{(-1)^r \rho^{(\beta)}(\xi(x,\bsy))\, \prod_{\ell = 1}^r D^{\bsalpha_\ell} \phi(\xi(x,\bsy), \bsy)}
{[D^0 \phi (\xi(x,\bsy), \bsy)]^{r + q}},
\end{equation}
which is of the form \cite[Formula~(3.10)]{GKSS:GKS22b} but now $(r, \bsalpha,
\beta)$ must satisfy $\bsmu + \beta \bse_0 + \sum_{\ell = 1}^r
\bsalpha_\ell = (r + q - 1, \bseta)$, noting the extra dependence on
$\bsmu$ due to the last condition in~\eqref{GKSS:eq:h-form}. Fortunately
the proof of \cite[Lemma~B.1]{GKSS:GKS22b} can be easily adapted to show that
the derivative with respect to $y_i$ of any function of the form
\eqref{GKSS:eq:h-tilde} can be written as
\begin{equation}
\label{GKSS:eq:D-h-tilde}
D^i \widetilde{h}_{q, \bseta, (r, \bsalpha, \beta,\bsmu)}(x,\bsy)
= \sum_{k = 1}^{\widetilde{K}_{q, \bseta}}
\widetilde{h}^{[k]}_{q, \bseta + \bse_i, (\widetilde{r}, \widetilde{\bsalpha}, \widetilde{\beta},\bsmu)}(x,\bsy),
\end{equation}
where $\widetilde{K}_{q, \bseta} = 4r + 2q + 1$ and each
$\widetilde{h}_{q, \bseta + \bse_i, (\widetilde{r}, \widetilde{\bsalpha},
\widetilde{\beta},\bsmu)} $ is again a function of the form
\eqref{GKSS:eq:h-tilde}, where now $(\widetilde{r}, \widetilde{\bsalpha},
\widetilde{\beta})$ satisfy $\bsmu + \widetilde{\beta} \bse_0 + \sum_{\ell
= 1}^{\widetilde{r}} \widetilde{\bsalpha}_\ell = (\widetilde{r} + q - 1,
\bseta + \bse_i)$.

Differentiating \eqref{GKSS:eq:h=Dtheta-h} with respect to $y_i$ using the
product and chain rules gives
\begin{align}
\label{GKSS:eq:Dh}
&D^i h_{q, \bseta, (r, \bsalpha, \beta, \bsmu)}(x,\bsy) \nonumber\\
&= \big[D^{\bsmu + \bse_i}\theta(\xi(x,\bsy), \bsy)
+ D^{\bsmu + \bse_0}\theta(\xi(x,\bsy), \bsy) \, D^i \xi(x,\bsy)\big]\,
\widetilde{h}_{q, \bseta, (r, \bsalpha, \beta, \bsmu)}(x,\bsy)
\nonumber\\
&\qquad + D^\bsmu \theta(\xi(x,\bsy), \bsy) \,
D^i \widetilde{h}_{q, \bseta, (r, \bsalpha, \beta, \bsmu)}(x,\bsy)
\nonumber\\
&= D^{\bsmu + \bse_i}\theta(\xi(x,\bsy), \bsy) \, \widetilde{h}_{q, \bseta,(r, \bsalpha, \beta, \bsmu)}(x,\bsy) \nonumber \\
&\qquad - \frac{ D^{\bsmu + \bse_0}\theta(\xi(x,\bsy), \bsy)  \,D^i \phi(\xi(x,\bsy), \bsy)}{D^0 \phi(\xi(x,\bsy), \bsy)}\,
\widetilde{h}_{q, \bseta,(r, \bsalpha, \beta, \bsmu)}(x,\bsy)
\nonumber\\
&\qquad + D^\bsmu \theta(\xi(x,\bsy), \bsy) \,
\sum_{k = 1}^{\widetilde{K}_{q, \bseta}}
\widetilde{h}^{[k]}_{q, \bseta + \bse_i, (\widetilde{r}, \widetilde{\bsalpha}, \widetilde{\beta},\bsmu)}(x,\bsy),
\end{align}
where we used \eqref{GKSS:eq:Di-xi} and \eqref{GKSS:eq:D-h-tilde}.

The first term in \eqref{GKSS:eq:Dh} can be expressed as $h_{q, \bseta +
\bse_i, (r, \bsalpha, \beta, \bsmu + \bse_i)}(x,\bsy)$.

The second term in \eqref{GKSS:eq:Dh} can be expressed as $h_{q, \bseta +
\bse_i, (r + 1, \widehat{\bsalpha}, \beta, \bsmu + \bse_0)}(x,\bsy)$,
where $\widehat{\bsalpha} = (\widehat{\bsalpha}_\ell)_{\ell = 1}^{r + 1}$
with $\widehat{\bsalpha}_\ell = \bsalpha_\ell$ for $\ell = 1, 2, \ldots,
r$ and $\widehat{\bsalpha}_{r + 1} = \bse_i$.

The third term in \eqref{GKSS:eq:Dh} is itself a sum of $\widetilde{K}_{q,
\bseta} = 4r+2q+1$ terms, where each term can be expressed in the form
$h_{q, \bseta + \bse_i, (\widetilde{r}, \widetilde{\bsalpha},
\widetilde{\beta}, \bsmu)}(x,\bsy)$.

Hence we conclude that \eqref{GKSS:eq:Dh} can be written as a sum of the
form \eqref{GKSS:eq:Dh-sum} with $K_{q, \bseta} \coloneqq \widetilde{K}_{q,
\bseta} + 2 = 4r + 2q + 3$ functions, where each function is of the form
\eqref{GKSS:eq:h-form} with $\bseta$ replaced by $\bseta + \bse_i$. The
upper bound $K_{q, \bseta} \leq 8|\bseta| + 6q - 1$ follows by the
restriction on $r$ in \eqref{GKSS:eq:h-form}.
\end{proof}

The following theorem generalises \cite[Theorem~3.2]{GKSS:GKS22b} to include
the extra function $\theta$ in \eqref{GKSS:eq:g_ind}. It can also be
viewed as an extension of \cite[Theorem~3]{GKSS:GKLS18} with a new bound on the
norm. For completeness we also incorporate the statement of
\cite[Theorem~3.3]{GKSS:GKS22b}.

\begin{theorem}\label{GKSS:thm:main2}
Let $d\ge 2$, $\bsnu \in \bbN_0^d$, and $x\in [a,b]\subset\R$. Suppose
that $\theta,\phi,\rho$ satisfy Assumption~\ref{GKSS:asm:phi} and
Assumption~\ref{GKSS:asm:theta-phi} for $q\in\{0,1\}$ and all
$\bseta\le\bsnu$.
\begin{itemize}
\item The function $g(y_0,\bsy) \coloneqq \theta(y_0,\bsy)\,\GKSSind(x
    - \phi(y_0,\bsy))$ satisfies $P_0 g \in \calH^{\bsnu}_{d} \cap
    C^{\bsnu}(\R^{d})$, with its $\calH^{\bsnu}_{d}$-norm bounded
    uniformly in $x$,
\begin{equation}
\label{GKSS:eq:P0-theta_norm}
 \|P_0 g\|_{\calH^\bsnu_{d}}
 \le
 \bigg(2\,\|\theta\|_{\calH^{(0,\bsnu)}_{d+1}}^2
 + 2 \sum_{\bszero\ne\bseta \leq \bsnu} \frac{\big(8^{|\bseta|-1}(|\bseta|-1)!\big)^2 B_{0,\bseta}}
 {\gamma_\bseta} \bigg)^{1/2}
 < \infty.
\end{equation}
\item In the case where $\theta \equiv 1$, we have
\begin{equation}
\label{GKSS:eq:P0-ind_norm}
\|P_0 g \|_{\calH^\bsnu_{d}} \le
\bigg( 1 + \sum_{\bszero \neq \bseta \leq \bsnu} \frac{\big(8^{|\bseta| - 1}(|\bseta| - 1)!\big)^2
B_{0, \bseta}}{\gamma_\bseta}\bigg)^{1/2}
< \infty.
\end{equation}
\item The function $g(y_0,\bsy) \coloneqq \delta(x - \phi(y_0,\bsy))$
    satisfies $P_0 g \in \calH^{\bsnu}_{d} \cap C^{\bsnu}(\R^{d})$,
    with
\begin{equation}
\label{GKSS:eq:P0-delta_norm}
 \|P_0 g\|_{\calH^\bsnu_{d}}
 \le \bigg(\sum_{\bseta \leq \bsnu} \frac{\big(8^{|\bseta|}|\bseta|!\big)^2 B_{1,\bseta}}
 {\gamma_\bseta} \bigg)^{1/2}
 < \infty.
\end{equation}
\end{itemize}
\end{theorem}

\begin{proof}
The second and third cases are precisely \cite[Theorem~3.2]{GKSS:GKS22b} and
\cite[Theorem~3.3]{GKSS:GKS22b}. Also, for $\theta \not \equiv 1$, the result
that $P_0 g \in \calH_d^\bsnu \cap C^\bsnu(\R^d)$ is given by
\cite[Theorem~3]{GKSS:GKLS18}. Hence it remains to prove the bound
\eqref{GKSS:eq:P0-theta_norm}. We focus on the non-trivial scenario when
$U_x$ defined in \eqref{GKSS:eq:Ux} is not empty.

Let $\bsy \in U_x$. For any $\bseta\in\bbN_0^d$ with
$\bszero\ne\bseta\le\bsnu$, we first prove by induction on $|\bseta|\ge 1$
that the $\bseta$th derivative of $P_0 g$ for the first case in
\eqref{GKSS:eq:P0g-simp} is given by
\begin{align} \label{GKSS:eq:hyp}
 D^\bseta [P_0 g(\bsy)]
 &=
 \int_{-\infty}^{\xi(x,\bsy)} D^\bseta\theta(y_0,\bsy)\,\rho_0(y_0)\rd y_0
 + \sum_{j=1}^{J_{0,\bseta}} h_{0,\bseta}^{[j]}(x,\bsy),
\end{align}
with $J_{0,\bseta} \le 8^{|\bseta|-1}(|\bseta|-1)!$ and where each
function $h_{0,\bseta}^{[j]}$ is of the form \eqref{GKSS:eq:h-form} with
$q=0$. For $\bsy \in \R^d \setminus U_x$, we have $\partial^\bseta
P_0g(\bsy) = 0$ for all $\bseta$.

For the base case $\bseta = \bse_i$ with any $i\in \{1:d\}$, using the
fundamental theorem of calculus and \eqref{GKSS:eq:Di-xi}, we obtain
\begin{align*}
  &D^i [P_0g(\bsy)]
  =
  \int_{-\infty}^{\xi(x,\bsy)}
  D^i\theta(y_0,\bsy)\,\rho_0(y_0) \rd y_0
  + \theta(\xi(x,\bsy),\bsy) \,\rho_0(\xi(x,\bsy))\, D^i\xi(x,\bsy) \nonumber \\
  &=
  \int_{-\infty}^{\xi(x,\bsy)} D^i\theta(y_0,\bsy)\,\rho_0(y_0) \rd y_0
  - \frac{\theta(\xi(x,\bsy),\bsy)\,\rho_0(\xi(x,\bsy))\,D^i \phi(\xi(x,\bsy),\bsy)}{D^0\phi(\xi(x,\bsy),\bsy)},
\end{align*}
which is recovered by taking $r=1$, $\bsalpha_1 = \bse_i$, $\beta=0$,
$\bsmu = \bszero$ and $J_{0,\bse_i}=1$ in \eqref{GKSS:eq:hyp}.

Suppose next that \eqref{GKSS:eq:hyp} holds for some $\bseta\in\bbN_0^d$
with $|\bseta|\ge 1$, and consider any $i\in \{1:d\}$ and $\bsy\in U_x$.
We have
\begin{align*}
 D^i D^\bseta [P_0g(\bsy)] =
 \int_{-\infty}^{\xi(x,\bsy)} D^{\bseta+\bse_i} \theta(y_0,\bsy)\,\rho_0(y_0) \rd y_0
 + \mbox{Extra},
\end{align*}
where
\begin{align*}
 \mbox{Extra}
 &\coloneqq
 D^\bseta \theta(\xi(x,\bsy),\bsy) \,\rho_0(\xi(x,\bsy))\, D^i\xi(x,\bsy)
 + \sum_{j=1}^{J_{0,\bseta}} D^i h_{0,\bseta}^{[j]}(x,\bsy) \\
 &\,=
 - \frac{D^\bseta\theta(\xi(x,\bsy),\bsy)\,\rho_0(\xi(x,\bsy))\,D^i \phi(\xi(x,\bsy),\bsy)}{D^0\phi(\xi(x,\bsy),\bsy)}
 + \sum_{j=1}^{J_{0,\bseta}} \sum_{k=1}^{K_{0,\bseta}} h_{0,\bseta+\bse_i}^{[j,k]}(x,\bsy) \\
 &= \sum_{j'=1}^{J_{0,\bseta+\bse_i}} h_{0,\bseta+\bse_i}^{[j']}(x,\bsy).
\end{align*}
In the second equality we used \eqref{GKSS:eq:Di-xi} as well as
Lemma~\ref{GKSS:lem:count} with $q=0$. The latter states that each
function $D^i h_{0,\bseta}^{[j]}(x,\bsy)$ can be written as a sum of
$K_{0,\bseta}\le 8|\bseta|-1$ functions of the form \eqref{GKSS:eq:h-form}
with $\bseta$ replaced by $\bseta+\bse_i$. We enumerated these functions
with the notation $h_{0,\bseta+\bse_i}^{[j,k]}(x,\bsy)$ and then relabeled
all functions for different combinations of indices $j$ and $k$ with the
notation $h_{0,\bseta+\bse_i}^{[j']}(x,\bsy)$. Also the fraction term in
the second equality can be expressed in the form
$h_{0,\bseta+\bse_i,(r,\bsalpha,\beta,\bsmu)}$ with $r=1$, $\bsalpha_1 =
\bse_i$, $\beta = 0$ and $\bsmu=\bseta$. The total number of functions
satisfies
\begin{align*}
 J_{0,\bseta+\bse_i} = 1 + J_{0,\bseta} K_{0,\bseta}
 \le 1 + 8^{|\bseta|-1}(|\bseta|-1)!\, (8|\bseta|-1)
 \le 8^{|\bseta+\bse_i|-1}\,(|\bseta+\bse_i|-1)!\,,
\end{align*}
as required. This completes the induction proof for \eqref{GKSS:eq:hyp}.

From \eqref{GKSS:eq:hyp} we have
\begin{align*}
 &\|P_0g\|_{\calH_{d}^{\bsnu}}^2
 = \sum_{\bseta \leq \bsnu} \frac{1}{\gamma_\bseta}
 \int_{\R^d} |D^\bseta [P_0g(\bsy)]|^2\,
 \bspsi_{\bseta}(\bsy_\bseta)\,\bsrho_{-\bseta}(\bsy_{-\bseta}) \rd \bsy
 \\
 &= \int_{U_x} \bigg|\int_{-\infty}^{\xi(x,\bsy)}
  \theta(y_0,\bsy) \rho(y_0) \rd y_0\,\bigg|^2\bsrho(\bsy)  \rd \bsy \\
 &\;\; +
 \sum_{\bszero\ne\bseta \leq \bsnu} \frac{1}{\gamma_\bseta}
 \int_{U_x} \bigg|\int_{-\infty}^{\xi(x,\bsy)}\!\!\!\!\!
  D^\bseta\theta(y_0,\bsy)\rho(y_0) \rd y_0
  + \sum_{j=1}^{J_{0,\bseta}} h_{0,\bseta}^{[j]}(\bsy)\bigg|^2
 \bspsi_{\bseta}(\bsy_\bseta)\,\bsrho_{-\bseta}(\bsy_{-\bseta}) \rd \bsy,
\end{align*}
and we estimate it by
\begin{align*}
 \|P_0g\|_{\calH_{d}^{\bsnu}}^2
 &\le \int_{\bbR^{d+1}} |\theta(y_0,\bsy)|^2 \rho(y_0) \bsrho(\bsy)  \rd (y_0,\bsy) \\
 &\quad + 2
 \sum_{\bszero\ne\bseta \leq \bsnu} \frac{1}{\gamma_\bseta}
  \int_{\bbR^{d+1}} |D^\bseta\theta(y_0,\bsy)|^2 \rho(y_0)\,
 \bspsi_{\bseta}(\bsy_\bseta)\,\bsrho_{-\bseta}(\bsy_{-\bseta}) \rd (y_0,\bsy)
 \\
 &\quad + 2
 \sum_{\bszero\ne\bseta \leq \bsnu} \frac{J_{0,\bseta}}{\gamma_\bseta}
 \sum_{j=1}^{J_{0,\bseta}} \int_{U_x} | h_{0,\bseta}^{[j]}(x,\bsy)|^2\,
 \bspsi_{\bseta}(\bsy_\bseta)\,\bsrho_{-\bseta}(\bsy_{-\bseta}) \rd \bsy
 \\
&\le 2\,\|\theta\|_{\calH^{(0,\bsnu)}_{d+1}}^2
 + 2 \sum_{\bszero\ne\bseta \leq \bsnu} \frac{\big(8^{|\bseta|-1}(|\bseta|-1)!\big)^2 B_{0,\bseta}}
 {\gamma_\bseta},
\end{align*}
where we used the Cauchy--Schwarz inequality and the simple estimate 
$|a + \sum_{k=1}^L b_k|^2 \le 2|a|^2 + 2|\sum_{k=1}^L b_k|^2 \le 2|a|^2 +
2L\sum_{k=1}^L |b_k|^2$, completing the proof.
\end{proof}

\section{QMC analysis for Asian option with preintegration}
\label{GKSS:sec:qmc-preint}

In this section we analyse the RMS error that results from approximating
the integrals \eqref{GKSS:eq:price}--\eqref{GKSS:eq:pdf} by a randomly
shifted lattice rule with a preintegration step. The lattice rules are
constructed by using the CBC construction, where the input weights are
specifically tailored to the $d$-dimensional integrands that are obtained
\emph{after} preintegration.

We now specialise to the function $\phi$ given by the average asset price
\eqref{GKSS:eq:phi}, and we set $\theta= x - \phi$ (and $x=K$) for
estimating the fair price of the Asian put option or $\theta\equiv 1$ for
estimating the cdf. In this case,  
we also have an explicit expression for the preintegrated fair price using \eqref{GKSS:eq:P_0} with $\theta=x-\phi$, 
\begin{align}
\label{GKSS:eq:P0g-price}
P_0&g(\bsy) 
\,=\, x\,\Phi(\xi(x,\bsy)) -\frac{S_0}{d+1}
\\\nonumber
&\cdot \bigg[\sum_{k=0}^{d} \exp \bigg( \big(R - \tfrac{1}{2}\sigma^2\big)\frac{(k+1)T}{d+1} + \sigma \bsA_{k,1:d}\bsy +\frac{\sigma^2A_{k,0}^2}{2} \bigg)\Phi(\xi(x,\bsy)-\sigma A_{k,0})\bigg],
\end{align}
where $\bsA_{k,1:d}$ is defined to be the vector $\bsA_k$ excluding the $A_{k,0}$ component.

Each integral \eqref{GKSS:eq:price}--\eqref{GKSS:eq:pdf} is now approximated
by $Q_{d, N} (P_0g)$ using one of the formulas \eqref{GKSS:eq:P0g-price} or 
\eqref{GKSS:eq:P0g-simp} for $P_0 g$, which are exact in the sense that
they no longer involve integrals, apart from the Gaussian cdf $\Phi$ 
(which can be approximated efficiently to machine precision). 
Then, to evaluate $P_0g(\bsPhi^{-1}(\bstau_k))$
for each QMC point $\bstau_k$, we compute the point of discontinuity
$\xi(x, \bsPhi^{-1}(\bstau_k))$ numerically using Newton's method.

\subsection{Bounds on the derivatives of the preintegrated functions}

In order to bound the norm of the preintegrated function $P_0 g$ in
\eqref{GKSS:eq:P0g-simp}, we need to know that
Assumption~\ref{GKSS:asm:phi} and Assumption~\ref{GKSS:asm:theta-phi} are
satisfied, and in addition we need to bound the integrals
\eqref{GKSS:eq:h-int} to obtain explicit values $B_{q, \bseta}$.

\begin{lemma} \label{GKSS:lem:B_pca}
Let $d\ge 1$, $\bsnu\in\bbN_0^d$ and $x\in [a,b]\subset\R$. Let $\phi$ be
the average asset price \eqref{GKSS:eq:phi} with $A$ given by the PCA
factorisation as in \eqref{GKSS:eq:A_pca}, take $\rho$ as the standard
normal density, let $\psi_i$ satisfy \eqref{GKSS:eq:psi-cond}, and set $\theta = 
x - \phi$ or $\theta\equiv 1$. Then
\begin{itemize}
\item Assumption~\ref{GKSS:asm:phi} holds; and
\item Assumption~\ref{GKSS:asm:theta-phi} holds for all
    $q\in\{0,1\}$ and $\bseta\le\bsnu$ satisfying $|\bseta|+q\le
    |\bsnu|+1$, with
\begin{align}
\label{GKSS:eq:B_pca}
 &B_{q, \bseta}
 \coloneqq \bigg(\frac{(\max_{\beta\le |\bseta|+q-1} \kappa_\beta)\,\Omega_q(2d+3)^{2|\bseta|+2q-1}}{
 [\min(\sigma \tau_d(2d+3),1)]^{|\bseta|+q}} \prod_{i = 1}^d
 \Lambda_i^{\eta_i} \bigg)^2, \\
 &\Lambda_i \coloneqq \frac{\sigma \tau_d (2d + 3)}{2i + 1}, \quad
 \Omega_q
 \coloneqq \begin{cases}
 1 & \mbox{if $\theta\equiv 1$ and $q=0$, or if $\theta= x- \phi$ and $q=1$}, \\
 \frac{1}{a} & \mbox{if $\theta\equiv 1$ and $q=1$}, \\
 b & \mbox{if $\theta = \phi - x$ and $q=0$},
 \end{cases} \nonumber
\end{align}
where $\kappa_\beta$ is defined in \eqref{GKSS:eq:rho-beta} below.
Note that $\tau_d (2d+3) = \sqrt{T(2 + \frac{1}{d+1})}\in
(\sqrt{2T},\sqrt{2.5T}\,]$ is bounded independently of $d$.
\end{itemize}
\end{lemma}

\begin{proof}
Since the formula \eqref{GKSS:eq:phi} for the time-discretised average
price $\phi$ is a finite sum of exponentials, clearly $\phi$ is analytic
on $\R^{d+1}$. Thus the last two items of Assumption~\ref{GKSS:asm:phi}
hold trivially. For $\bseta \in \N_0^{d+1}$ and writing $(y_0,\bsy) \in
\R^{d+1}$, by differentiating \eqref{GKSS:eq:phi} we obtain
\begin{equation}
\label{GKSS:eq:d-phi}
D^\bseta \phi(y_0,\bsy)
= \frac{1}{d+1} \sum_{k = 0}^d \bigg( \prod_{i = 0}^d [\sigma A_{k, i}]^{\eta_i}\bigg)
 S_0 \exp \bigg( \big(R - \tfrac{1}{2}\sigma^2\big)\frac{(k+1)T}{d+1} + \sigma \bsA_{k} (y_0,\bsy)\bigg).
\end{equation}
Since $\sigma > 0$ and the entries of $A$ given by \eqref{GKSS:eq:A_pca}
satisfy $A_{k, 0} > 0$ for all $k = 0, 1,\ldots d$, it follows from
\eqref{GKSS:eq:d-phi} with $\bseta = \bse_0$ that $D^0\phi(y_0,\bsy) > 0$
for all $(y_0,\bsy) \in \R^{d+1}$. Similarly, for each $\bsy\in \R^d$,
$\sigma
> 0$ and $A_{k, 0} > 0$ imply that $\phi(y_0,\bsy) \to \infty$ as
$y_0\to \infty$, because we can simply separate out $\exp(\sigma A_{k , 0}
y_0)$ in each term in the sum \eqref{GKSS:eq:phi}, which clearly goes to
$\infty$ as $y_0 \to \infty$. Hence, $\phi$ as in \eqref{GKSS:eq:phi}
satisfies the first two items of Assumption~\ref{GKSS:asm:phi}.

To verify Assumption~\ref{GKSS:asm:theta-phi}, we need to consider
functions of the form \eqref{GKSS:eq:h-form}, and in particular,
\begin{align} \label{GKSS:eq:h-abs}
 |h_{q,\bseta}
 (x,\bsy)|
 =
\frac{|D^\bsmu\theta(\xi(x,\bsy), \bsy)| \ |\rho^{(\beta)}(\xi(x,\bsy))|
\prod_{\ell = 1}^r |D^{\bsalpha_\ell} \phi(\xi(x,\bsy), \bsy)|} {|D^0
\phi(\xi(x,\bsy), \bsy)|^{r+q}}.
\end{align}
Since $\rho$ is the standard normal density, the order $\beta$ derivative
is $\rho^{(\beta)}(u)  = (-1)^\beta \GKSSHe_\beta(u) \rho(u)$, where
$\GKSSHe_\beta$ is the order $\beta$ probabilist's Hermite polynomial
which is in turn dominated by the density $\rho$, so that
$\rho^{(\beta)}(u) \to 0$ as $u \to \pm\infty$ and 
$\kappa_\beta < \infty$ is defined such that, for all $u \in \R$
\begin{equation}
\label{GKSS:eq:rho-beta}
|\rho^{(\beta)}(u)| \,=\, \GKSSHe_\beta(u) \rho(u) \,\leq\, 
\kappa_\beta \,\coloneqq\, \sup_{v \in \R} \GKSSHe_\beta(v) \rho(v).
\end{equation}
Similarly, all derivatives of $\phi$ and $\theta$ will be
dominated by $\rho$. Thus \eqref{GKSS:eq:h-lim} and \eqref{GKSS:eq:h-int}
hold.

It remains to obtain explicit values of $B_{q,\bseta}$ in
\eqref{GKSS:eq:h-int}. We begin by obtaining an upper bound to
$|D^{\bsalpha_\ell}\phi(\xi(x,\bsy),\bsy)|$ in the numerator of
\eqref{GKSS:eq:h-abs}. Using $\sin(\tfrac{\pi}{2} u) \geq u$ for $u
\in [0,1]$ and $\sin(v) \leq 1$ for $v \in \R$, the matrix entries
in \eqref{GKSS:eq:A_pca} satisfy 
\begin{align} \label{GKSS:eq:A_ki_upper}
\sigma |A_{k, i}|
\le \sigma \tau_d\, \bigg|\frac{\sin(2(k+1)\chi_{i})}{\sin(\chi_{i})}\bigg|
\le \frac{\sigma \tau_d (2d + 3)}{2i + 1} \eqqcolon \Lambda_i.
\end{align}
It then follows from \eqref{GKSS:eq:d-phi} with $\bseta = \bsalpha_\ell$
and the definition of $\xi$ in \eqref{GKSS:eq:xi} that 
\begin{align}
\label{GKSS:eq:dphi_pca_upper}
 |D^{\bsalpha_\ell} \phi(\xi(x,\bsy),\bsy)|
 &\le \big(\textstyle\prod_{i = 0}^d \Lambda_i^{\alpha_{\ell,i}}\big)\,
 \phi(\xi(x,\bsy),\bsy)
 = \big(\textstyle\prod_{i = 0}^d \Lambda_i^{\alpha_{\ell,i}}\big)\,x.
\end{align}

Next we obtain a lower bound to $|D^0 \phi(\xi(x,\bsy),\bsy)|$ in the
denominator of \eqref{GKSS:eq:h-abs}. Using $\chi_0 = \frac{\pi}{2(2d+3)}$
from \eqref{GKSS:eq:A_pca} and $u \leq \sin(\tfrac{\pi}{2} u) \leq 2u$ for
$u\in [0,1]$, we obtain
\begin{align*}
 \sigma A_{k,0}
 = \sigma\tau_d\, \frac{\sin(2(k+1)\chi_0)}{\sin(\chi_0)}
 \ge \sigma\tau_d  (k+1)
 \ge \sigma\tau_d.
\end{align*}
Thus it follows from \eqref{GKSS:eq:d-phi} with $\bseta = \bse_0$ and
again \eqref{GKSS:eq:xi} that
\begin{equation}
\label{GKSS:eq:d1phi_pca_lower}
D^0 \phi(\xi(x,\bsy),\bsy) \ge \sigma \tau_d\,\phi(\xi(x,\bsy),\bsy)
= \sigma \tau_d\, x.
\end{equation}

Finally we consider $|D^\bsmu\theta(\xi(x,\bsy),\bsy))|$ in the numerator
of \eqref{GKSS:eq:h-abs}. Recall that we have either $\theta\equiv 1$ or
$\theta = x - \phi$, and
\begin{align} \label{GKSS:eq:Dmu}
 |D^\bsmu\theta(\xi(x,\bsy),\bsy))|
 &= \begin{cases}
 1 & \mbox{if $\theta\equiv 1$ and $\bsmu=\bszero$}, \\
 0 & \mbox{if $\theta\equiv 1$ and $\bsmu\ne\bszero$}, \\
 x - \phi(\xi(x,\bsy), \bsy)=0 & \mbox{if $\theta= x - \phi$ and $\bsmu=\bszero$}, \\
 |D^\bsmu\phi(\xi(x,\bsy), \bsy)|
 \le
 \big(\prod_{i = 0}^d \Lambda_i^{\mu_i}\big)\,x\quad
  & \mbox{if $\theta= x - \phi$ and $\bsmu\ne\bszero$}.
\end{cases}
\end{align}
Observe that all four cases can be bounded by the last case, but with $x$
replaced by $x^m$, where $m\coloneqq 0$ if $\theta\equiv 1$ and
$m \coloneqq 1$ if $\theta =  x - \phi$. (Note that $\bsmu=\bszero$ in
the first case means the product factor is exactly $1$.)

Combining \eqref{GKSS:eq:rho-beta}, \eqref{GKSS:eq:dphi_pca_upper},
\eqref{GKSS:eq:d1phi_pca_lower} and \eqref{GKSS:eq:Dmu}, we obtain
\begin{align*}
 &|h_{q,\bseta}(x,\bsy)|
 \le \frac{\kappa_\beta}{[\sigma \tau_d\, x]^{r+q}}
 \bigg(\prod_{i = 0}^d \Lambda_i^{\mu_i}\bigg)\,x^m
 \prod_{\ell=1}^r \bigg[\bigg(\prod_{i = 0}^d \Lambda_i^{\alpha_{\ell,i}}\bigg) x\bigg]\\
 &= \frac{\kappa_\beta x^{m-q}[\sigma\tau_d(2d+3)]^{\mu_0 + \sum_{\ell=1}^r \!\!\alpha_{\ell,0}}}{[\sigma \tau_d]^{r+q}}
 \prod_{i = 1}^d \Lambda_i^{\mu_i + \sum_{\ell=1}^r \!\!\alpha_{\ell,i}}
 = \frac{\kappa_\beta x^{m-q}(2d+3)^{r+q}}{[\sigma \tau_d(2d+3)]^{1+\beta}}
 \prod_{i = 1}^d \Lambda_i^{\eta_i},
\end{align*}
where we used $\bsmu + \beta\bse_0 + \sum_{\ell=1}^r \bsalpha_\ell =
(r+q-1, \bseta)$. Next we seek an upper bound that is independent of
$(r,\bsalpha,\beta, \bsmu)$. Using $r\le 2|\bseta|+q-1$ and $\beta\le
|\bseta|+q-1$, and noting that 
$u^{1 + \beta} \geq \min(u, 1)^{1 + \beta} \geq \min(u, 1)^{|\bseta| + q}$ for $ u > 0$,
we arrive at
\begin{align*}
 |h_{q,\bseta}(x,\bsy)|
 \le \frac{(\max_{\beta\le |\bseta|+q-1} \kappa_\beta)\,
 x^{m-q}(2d+3)^{2|\bseta|+2q-1}}{
 [\min(\sigma \tau_d(2d+3),1)]^{|\bseta|+q}}
 \prod_{i = 1}^d \Lambda_i^{\eta_i}.
\end{align*}
Since the above upper bound on $|h_{q,\bseta}(x,\bsy)|$ is independent of
$\bsy$, and $\rho$ and all $\psi_i$ integrate to $1$, it is easy to see
that the integral \eqref{GKSS:eq:h-int} is bounded by $B_{q,\bseta}$ as
given in \eqref{GKSS:eq:B_pca}, where $\Omega_q$ is defined such that
$x^{m-q}\le \Omega_q$ with $x\in [a,b]$.
\end{proof}

\begin{remark}
Note that for $i \neq 0$, the entries $A_{k, i}$ given by \eqref{GKSS:eq:A_pca} 
are both positive and negative, and so $\phi$
is not monotone with respect to any variable except $y_0$.
Hence, for the PCA factorisation $y_0$ is the only possible choice for the preintegration variable.
For Cholesky and Brownian Bridge, $\phi$ is monotone with respect to
every variable and so a similar analysis of preintegration could be 
performed in those cases as well. 
\end{remark}

\subsection{Choosing the weights}

We showed in the previous subsection that the fair price, cdf and pdf
associated with Asian option satisfy the conditions for preintegration,
and hence a randomly shifted lattice rule can be constructed such that the
approximation $Q_{d, N}(P_0 g)$ achieves nearly order $1/N$ convergence
rate. Now we tailor the lattice rule to each integrand $P_0 g$ by choosing
weights which minimise the resulting bound on the RMS error
\eqref{GKSS:eq:cbc-err}.

\begin{corollary}
\label{GKSS:cor:norm}
Under the setting of Lemma~\ref{GKSS:lem:B_pca}, the preintegrated
functions $P_0 g$ \eqref{GKSS:eq:P0g-simp} corresponding to the three cases in
Theorem~\ref{GKSS:thm:main2} satisfy
\begin{align*}
 \|P_0 g\|_{\calH^{\bsnu}_d} \le \bigg(\sum_{\bseta\le\bsnu}
 \frac{A_\bseta}{\gamma_\bseta}\bigg)^{1/2},
 \;\mbox{with}\;
 A_\bszero \coloneqq
 \begin{cases}
 2 Z^2 \prod_{i=0}^d I_{1,i} & \mbox{if $g = (x-\phi)\ind(x-\phi)$,} \\
 1 & \mbox{if $g = \ind(x-\phi)$,} \\
 B_{1,\bszero} & \mbox{if $g = \delta(x-\phi)$,} \\
 \end{cases}
\end{align*}
and for $\bseta\ne\bszero$,
\begin{align*}
 A_\bseta \coloneqq
 \begin{cases}
 2 Z^2 I_{1,0} \big(\prod_{i = 1, \eta_i= 0}^d I_{1,i}\big)
 \prod_{i = 1, \eta_i\ne 0}^d \big(\Lambda_i^{2\eta_i} I_{2,i}\big)
 + 2\!\!\!\!\!\!\! &
 \big(8^{|\bseta|-1}(|\bseta|-1)!\big)^2 B_{0,\bseta} \\
 & \mbox{if $g = (x-\phi)\ind(x-\phi)$,} \\
 \big(8^{|\bseta|-1}(|\bseta|-1)!\big)^2 B_{0,\bseta} & \mbox{if $g = \ind(x-\phi)$,} \\
 \big(8^{|\bseta|}|\bseta|!\big)^2 B_{1,\bseta} & \mbox{if $g = \delta(x-\phi)$,} \\
 \end{cases}
\end{align*}
where $Z$ is defined in \eqref{GKSS:eq:defZ} below,
\begin{align*}
 I_{1,i} \coloneqq
 \int_{-\infty}^\infty
 e^{2\Lambda_i|y|}\rho(y)\rd y = 2e^{2\Lambda_i^2} \Phi(2\Lambda_i)
 \quad\mbox{and}\quad
 I_{2,i} \coloneqq
 \int_{-\infty}^\infty
 e^{2\Lambda_i|y|}\psi_i(y)\rd y.
\end{align*}
\end{corollary}

\begin{proof}
Lemma~\ref{GKSS:lem:B_pca} ensures that all the required assumptions are
satisfied for Theorem~\ref{GKSS:thm:main2}, and hence we have the three
bounds \eqref{GKSS:eq:P0-theta_norm}, \eqref{GKSS:eq:P0-ind_norm},
\eqref{GKSS:eq:P0-delta_norm}. In particular, for the bound
\eqref{GKSS:eq:P0-theta_norm} with $\theta= x - \phi$ we also need
\begin{align*}
 \|\theta\|_{\calH_{d+1}^{(0,\bsnu)}}^2
 &= \sum_{\bseta \leq \bsnu} \frac{1}{\gamma_\bseta}
 \int_{\R^{d+1}} |D^{(0,\bseta)} \phi(y_0,\bsy)]|^2\,
 \bspsi_{\bseta}(\bsy_\bseta)\,\bsrho_{-\bseta}(\bsy_{-\bseta})\,\rho(y_0)\rd y_0 \rd \bsy
 \\
 &\le \sum_{\bseta \leq \bsnu} \frac{1}{\gamma_\bseta}
 \bigg(\prod_{i = 1}^d \Lambda_i^{\eta_i}\bigg)^2
 \int_{\R^{d+1}}  | \phi(y_0,\bsy) |^2\,
 \bspsi_{\bseta}(\bsy_\bseta)\,\bsrho_{-\bseta}(\bsy_{-\bseta})\,\rho(y_0)\rd y_0 \rd \bsy, 
\end{align*}
where we adapted \eqref{GKSS:eq:dphi_pca_upper}. Using \eqref{GKSS:eq:phi}
and \eqref{GKSS:eq:A_ki_upper} we obtain
\begin{align} \label{GKSS:eq:defZ}
 | \phi(y_0,\bsy) |
 &\le \underbrace{\frac{1}{d+1} \sum_{k = 0}^d S_0
 \exp \bigg( \big(R - \tfrac{1}{2}\sigma^2\big)\frac{(k+1)T}{d+1}\bigg)}_{\eqqcolon\,Z}
 \prod_{i=0}^d e^{\Lambda_i |y_i|}.
\end{align}
Thus we conclude that
\begin{align*}
 \|\theta\|_{\calH_{d+1}^{(0,\bsnu)}}^2
 &\le \sum_{\bseta \leq \bsnu} \frac{Z^2\,I_{1,0}}{\gamma_\bseta}
 \bigg(\prod_{i = 1, \eta_i= 0}^d I_{1,i}\bigg)
 \prod_{i = 1, \eta_i\ne 0}^d \Big(\Lambda_i^{2\eta_i} I_{2,i}\Big),
\end{align*}
as required.
\end{proof}

To ensure that the integrals $I_{2,i}$ are finite, we choose  exponential weight functions
\begin{equation}
\label{GKSS:eq:psi-exp}
 \psi_i(y) \coloneqq \Lambda_0\, e^{- 2\Lambda_0 |y|}, \quad\mbox{with}\quad \Lambda_0 \coloneqq \sigma\tau_d(2d+3),
\end{equation}
so that $\int_{-\infty}^\infty \psi_i(y)\rd y = 1$. This leads to
$I_{2,i} = \Lambda_0 \int_{-\infty}^\infty e^{-2(\Lambda_0-\Lambda_i)|y|}
\rd y = 1 + \frac{1}{2i}$.

For this choice of $\psi_i$,
in the CBC error bound \eqref{GKSS:eq:cbc-err} with $h = P_0 g$
we can take $ \omega = 1 - \delta^2$  and $\lambda = 1/(2(1 - \delta))$ for $\delta \in (0, \tfrac{1}{2})$ (cf. \cite{GKSS:KSWW10}). We then 
choose the weights to minimise this bound using \cite[Lemma~6.2]{GKSS:KSS12}.

\begin{corollary}
Under the setting of Lemma~\ref{GKSS:lem:B_pca},
the error for the randomly shifted lattice rule constructed by the CBC algorithm \cite{GKSS:NK14} applied to each of  the 
three preintegrated functions \eqref{GKSS:eq:P0g-simp} satisfies
\begin{equation}
\label{GKSS:eq:err-final}
\mathrm{RMSE} \,\leq\, 
C_{\bsgamma, \delta} \, [\varphi_\mathrm{tot}(N)]^{-1 + \delta} \quad 
\text{for $\delta \in (0, \tfrac{1}{2})$},
\end{equation}
where $C_{\bsgamma, \delta}$ corresponds to \eqref{GKSS:eq:cbc-err}
with $h = P_0 g$ and the norm bounded as in Corollary~\ref{GKSS:cor:norm}.
The weights and constant minimising this upper bound are
\begin{align}
\label{GKSS:eq:gamma-opt}
\gamma_\bseta^* &\coloneqq \bigg(\frac{A_\bseta}
{[2C_2 \zeta(1 + \delta)]^{|\bseta|}} \bigg)^{\frac{2(1 - \delta)}{3 - 2\delta}},\\
\quad
C_{\bsgamma, \delta}^*& =
\bigg(\sum_{\bseta \in \{0, 1\}^d} \big[ A_\bseta^{\frac{1}{2(1 - \delta)}}
[2C_2 \zeta\big(1 + \delta)]^{|\bseta|}\big]^\frac{2(1 - \delta)}{3 - 2\delta}
\bigg)^{\frac{3}{2} - \delta} \notag,
\end{align}
i.e., $\{\gamma_\bseta\}$ are of ``product and order dependent'' (POD).
\end{corollary}

\section{Numerical results}
\label{GKSS:sec:num}
In this section we present convergence results for
approximating the fair price, the cdf and the pdf as in \eqref{GKSS:eq:price}--\eqref{GKSS:eq:pdf}.
Motivated by the choice of POD weights (41) with the three possible $A_{\bseta}$ given in Corollary \ref{GKSS:cor:norm}, here we take simplified product weights (corresponding to the limiting case of $\delta = 0$):
\begin{equation}
\label{GKSS:eq:weights-num}
 \gamma_\bseta = \bigg(\prod_{i=1}^d \Lambda_i^{2\eta_i}
 \bigg)^{2/3}
 = \prod_{i=1,\eta_i\ne 0}^d \Lambda_i^{4/3}.
\end{equation}
The error bound \eqref{GKSS:eq:err-final} is still valid, but the constant is not
minimised. 

We use the parameters
$K = \$100$, $S_0 = \$ 100$, $R = 0.1$, $\sigma = 0.2$, $T = 1$ and $d + 1 = 256$ time steps.
As a benchmark, we compare to plain MC and QMC, as well as MC with preintegration. 
For all QMC approximations we take the final estimate to be the average
over $L = 32$ random shifts, and for a fair comparison each MC
approximation uses $L\times N$ points in total. For all QMC approximations
we use the same generating vector given by the CBC construction \cite{GKSS:NK14}
with weights given by \eqref{GKSS:eq:weights-num} and weight functions
\eqref{GKSS:eq:psi-exp} for 
$N = 101, 251, 503, 997, 1999, 4001, 8009, 16001, 32003, 64007 ,128021$.
Since in each problem \eqref{GKSS:eq:price}--\eqref{GKSS:eq:pdf}
the original integrand $g$ is either discontinuous or has a kink,
there is no theory to construct a `good' generating vector for plain QMC.
As such, for plain QMC we use the same generating vector as for QMC after preintegration, 
except including an extra $(d+1)$th component.
The plain QMC approximations also use the PCA factorisation.
We estimate the RMSE by the sample standard error over the
random shifts, and for MC, the RMSE is estimated by the
sample standard error over all MC realisations.

\begin{figure}[!t]
\centering
\includegraphics[scale=0.43]{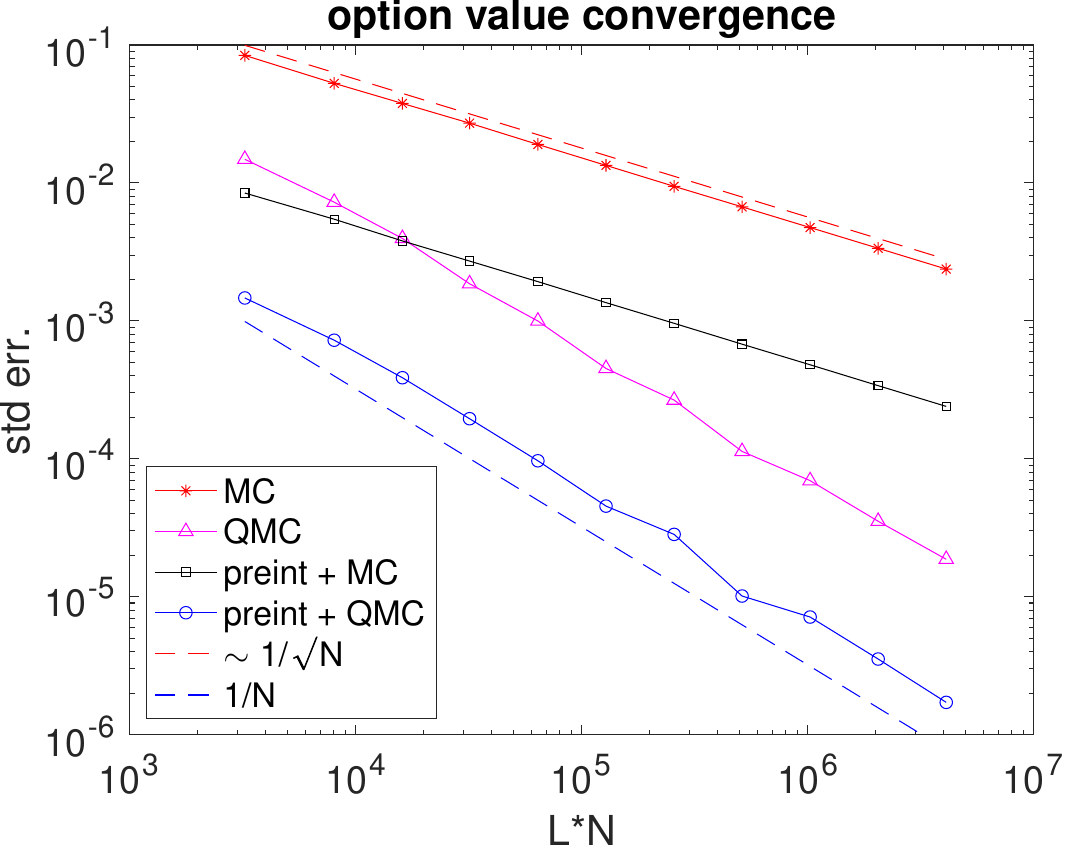}
\hfill
\includegraphics[scale=0.43]{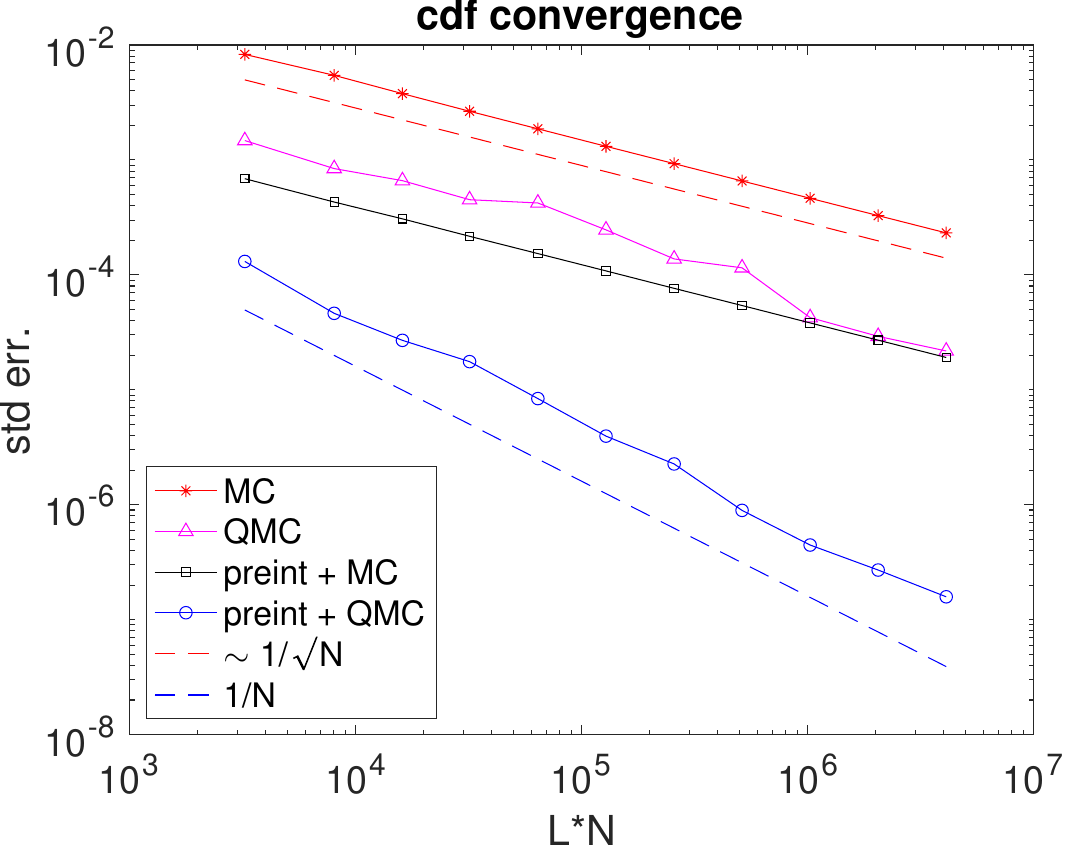}
\\[2mm]
\includegraphics[scale=0.43]{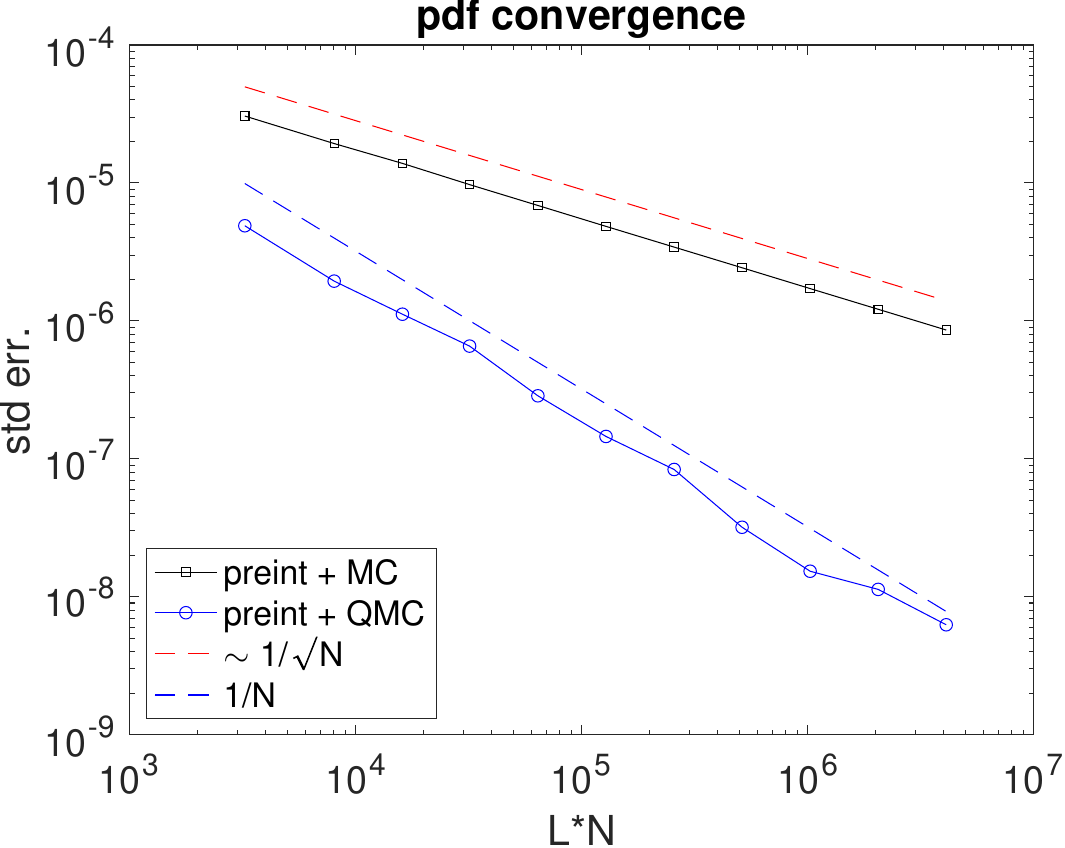}
\hfill
\includegraphics[scale=0.43]{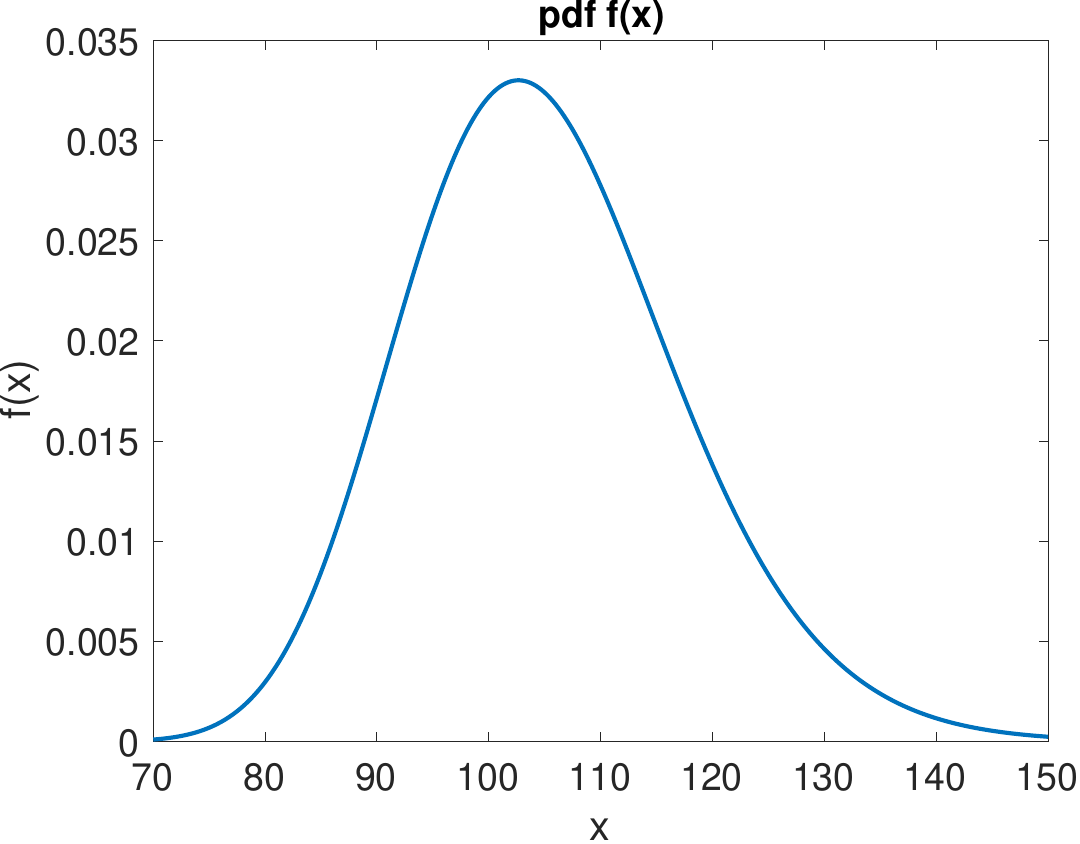}

\caption{RMSE convergence for MC, QMC, MC with preintegration and QMC with preintegration
for approximating: the fair price (top left), the cdf $F(x)$ for $x = 100$ (top right)
and the pdf $f(x)$ for $ x = 100$ (bottom left). Plot of approximated pdf $f(x)$ for $x \in [70, 150]$ (bottom right).
}
\label{GKSS:fig:opt:pointwise}
\end{figure}

In Figure~\ref{GKSS:fig:opt:pointwise} we study the convergence in $N$
of the estimated RMSE
for plain MC, plain QMC, MC with preintegration and QMC with preintegration.
In the top left we plot the RMSE for the fair price; 
in the top right we plot the RMSE for cdf at $x = K = 100$
(i.e., the probability of being ``in the money'');
and in the bottom left we plot the RMSE for the pdf at $x = 100$,
note that the Dirac $\delta$ cannot be evaluated in practice
and so plain MC and QMC approximations of the pdf are not possible.
In the bottom right of Figure~\ref{GKSS:fig:opt:pointwise}
we plot the approximated pdf  $f(x)$ for $x \in [70, 150]$, using $N = 16001$ QMC points 
and Chebyshev inperpolation with 30 points.

It is very clear from Figure~\ref{GKSS:fig:opt:pointwise} that for all three problems 
QMC with preintegration beats both
plain QMC and MC, both in absolute terms at all values of $N$, and also
with respect to rate of convergence.  We also note that plain QMC
manifestly beats plain MC, which is to be expected and
has been shown empirically many times previously
for similar option pricing problems.

In all of our experiments, the computational time required for QMC with preintegration
increased by less than a factor of 2 compared to plain QMC without preintegration,
whereas the error is reduced by a factor of ranging from 10 to 100.

\bibliographystyle{plain}

\end{document}